\documentclass[preprint,10pt]{elsarticle}
\pdfoutput=1
\usepackage{titlesec}



\makeatletter
\def\ps@pprintTitle{%
   \let\@oddhead\@empty
   \let\@evenhead\@empty
   \def\@oddfoot{\reset@font\hfil\thepage\hfil}
   \let\@evenfoot\@oddfoot
}
\makeatother

\usepackage{titlesec}
\usepackage{xcolor}
\usepackage{enumitem}
\usepackage{hyperref}
\hypersetup{colorlinks=true,urlcolor=black}
\usepackage[T1]{fontenc}\usepackage{lmodern}

\usepackage{mathtools,amsfonts,amssymb}
\usepackage{caption}
\usepackage{subcaption}
\usepackage{booktabs}
\usepackage{pifont}
\usepackage{placeins}

\usepackage{tabularx}

\definecolor{green}{rgb}{0,0.5,0}

\usepackage{tikz}
\usetikzlibrary{calc}
\usepackage{tkz-euclide}
\usepackage{calc}
\usepackage{tikz,pgfplots}
\pgfplotsset{compat=1.12}
\usetikzlibrary{arrows.meta,3d,matrix,positioning,decorations.markings}

\newcommand{\nmo}{\xi}
\newcommand{\ope}{\eta}
\newcommand{\npe}{\zeta}

\newcommand{\nbasis}{n}
\newcommand{\srcidx}{j}
\newcommand{\testidx}{i}

\makeatletter
\tikzoption{canvas is xy plane at z}[]{%
  \def\tikz@plane@origin{\pgfpointxyz{0}{0}{#1}}%
  \def\tikz@plane@x{\pgfpointxyz{1}{0}{#1}}%
  \def\tikz@plane@y{\pgfpointxyz{0}{1}{#1}}%
  \tikz@canvas@is@plane
}
\makeatother

\makeatletter
\tikzoption{canvas is plane}[]{\@setOxy#1}
\def\@setOxy O(#1,#2,#3)x(#4,#5,#6)y(#7,#8,#9)%
  {\def\tikz@plane@origin{\pgfpointxyz{#1}{#2}{#3}}%
   \def\tikz@plane@x{\pgfpointxyz{#4}{#5}{#6}}%
   \def\tikz@plane@y{\pgfpointxyz{#7}{#8}{#9}}%
   \tikz@canvas@is@plane
  }
\makeatother 

\definecolor{rev_red}{RGB}{227,26,28}

\definecolor{orange}{rgb}{1,0.5,0}
\definecolor{green}{rgb}{0,0.5,0}
\definecolor{purple}{rgb}{0.5,0,0.5}

\newcommand{\editor}[1]{#1}
\newcommand{\reviewerTwo}[1]{#1}
\usepackage[margin=1in]{geometry}
\biboptions{sort&compress}
\begin{document}

\begin{center}
{\Large
Code Verification for Practically Singular Equations}
\\[2em]
Brian A.~Freno$^\text{a}$ and Neil R.~Matula$^\text{a}$
\\[1em]
{\footnotesize\textit{%
$^\text{a}$Sandia National Laboratories, Albuquerque, NM 87185}}
\end{center}

%
%
%
%
%
%
%
%
%
%
%


\section{Introduction}

The method-of-moments implementation of the electric-field integral equation (EFIE) yields many code-verification challenges due to the various sources of numerical error and their possible interactions~\cite{freno_em_mms_2020,freno_em_mms_quad_2021}.  Matters are further complicated by singular integrals, which arise from the presence of a Green's function.  To address these singular integrals, an approach is presented in~\cite{freno_em_mms_2020} wherein both the solution and Green's function are manufactured.  Because the arising equations are poorly conditioned, they are reformulated as a set of constraints for an optimization problem that selects the solution closest to the manufactured solution.  
In this paper, we demonstrate how, for such practically singular systems of equations, computing the truncation error by inserting the exact solution into the discretized equations cannot detect certain orders of coding errors.  On the other hand, the discretization error from the optimal solution~\cite{freno_em_mms_2020} is a more sensitive metric \reviewerTwo{that can detect orders less than those of the expected convergence rate}.
\section{Truncation and Discretization Error}

Let $\boldsymbol{\mathcal{L}}(\cdot)$ and $\boldsymbol{\mathcal{L}}_h(\cdot)$ respectively denote the operators representing the continuous and discretized equations.  The truncation error $\boldsymbol{\tau}_h(\cdot)$ is the difference between them~\cite{ferziger_2002,oberkampf_2010}:
\begin{align}
\boldsymbol{\mathcal{L}}_h(\mathbf{v}) = \boldsymbol{\mathcal{L}}(\mathbf{v}) + \boldsymbol{\tau}_h(\mathbf{v}).
\label{eq:trunc}
\end{align}

Let $\mathbf{u}$ denote the solution to the continuous equations ($\boldsymbol{\mathcal{L}}(\mathbf{u})=\mathbf{0}$) and $\mathbf{u}_h$ denote the solution to the discretized equations ($\boldsymbol{\mathcal{L}}_h(\mathbf{u}_h)=\mathbf{0}$).
Inserting a Taylor-series expansion of $\boldsymbol{\mathcal{L}}(\mathbf{v})$ about $\mathbf{v}=\mathbf{u}$ into~\eqref{eq:trunc} and evaluating at $\mathbf{v}=\mathbf{u}_h$ yields
\begin{align}
\left.\frac{\partial \boldsymbol{\mathcal{L}}}{\partial\mathbf{v}}\right|_\mathbf{u}\!\cdot\mathbf{e}_h + \mathcal{O}(\|\mathbf{e}_h\|^2) =- \boldsymbol{\tau}_h(\mathbf{u}_h),
\label{eq:linearization}
\end{align}
where $\mathbf{e}_h=\mathbf{u}_h-\mathbf{u}$ is the discretization error.  Therefore, \eqref{eq:linearization} relates the truncation error to the discretization error.  

\reviewerTwo{%
The truncation error can be computed from~\eqref{eq:trunc} by setting $\mathbf{v}=\mathbf{u}$:
\begin{align}
\boldsymbol{\tau}_h(\mathbf{u}) = \boldsymbol{\mathcal{L}}_h(\mathbf{u}).
\label{eq:trunc1}
\end{align}
As described in Section 5.5.6.1 of~\cite{oberkampf_2010}, evaluating~\eqref{eq:trunc1} across a series of systematically refined discretizations can be used to compute the convergence rate of $\mathbf{e}_h$ without solving for $\mathbf{u}_h$.  However, there are known drawbacks to using this approach~\cite{oberkampf_2010}.}

\section{The Method-of-Moments Implementation of the EFIE}

The variational form of the EFIE is: find $\mathbf{u}$, such that
\begin{align}
a(\mathbf{u},\mathbf{v}) = \left(\mathbf{f}, \mathbf{v}\right) \qquad \forall \mathbf{v}\in\mathbb{V},
\label{eq:var_sesquilinear}
\end{align}
where $\mathbf{f}$ is a source, $\mathbb{V}$ is the \reviewerTwo{test} space, and the sesquilinear form and inner product are defined by
\begin{align}
a(\mathbf{u},\mathbf{v}) &{}=  \alpha \int_S \bar{\mathbf{v}}(\mathbf{x})\cdot\int_{S'} \mathbf{u}(\mathbf{x}')G(\mathbf{x},\mathbf{x}')dS'dS+ \beta \int_S \nabla\cdot\bar{\mathbf{v}}(\mathbf{x})\int_{S'} \nabla'\cdot\mathbf{u}(\mathbf{x}')G(\mathbf{x},\mathbf{x}')dS' dS, \label{eq:a}
\\ \nonumber
(\mathbf{u},\mathbf{v})  &{}= \int_S \mathbf{u}(\mathbf{x})\cdot \bar{\mathbf{v}}(\mathbf{x}) dS.
\end{align}
In~\eqref{eq:a}, $\alpha$ and $\beta$ are constants, $G$ is a Green's function, and the bar notation denotes complex conjugation.  If the domain is open, the component of $\mathbf{u}$ normal to the boundary is zero.  

To solve~\eqref{eq:var_sesquilinear}, we approximate $\mathbf{u}$ with $\mathbf{u}_h$ using a Galerkin method with Rao--Wilton--Glisson (RWG) basis functions $\boldsymbol{\phi}_{\srcidx}(\mathbf{x})$~\cite{rao_1982}:
%
$\mathbf{u}_h(\mathbf{x}) =\sum_{\srcidx=1}^{\nbasis} u_{\srcidx}^h \boldsymbol{\phi}_{\srcidx}(\mathbf{x})$,
%
where $\nbasis$ is the total number of unknowns.  Defining $\mathbb{V}_h$ to be the span of the RWG basis functions, the Galerkin approximation of the original problem is now: find $\mathbf{u}_h\in\mathbb{V}_h$, such that
\begin{align}
a(\mathbf{u}_h,\boldsymbol{\phi}_{\testidx}) = \left(\mathbf{f}, \boldsymbol{\phi}_{\testidx}\right)
\label{eq:proj_disc}
\end{align}
for $i=1,\hdots,\nbasis$.

To verify the order of accuracy of~\eqref{eq:proj_disc}, we can manufacture $\mathbf{f}$~\cite{freno_em_mms_2020} so that~\eqref{eq:proj_disc} becomes
\begin{align}
a(\mathbf{u}_h,\boldsymbol{\phi}_{\testidx}) = a(\mathbf{u},\boldsymbol{\phi}_{\testidx}).
\label{eq:proj_disc2}
\end{align}
%
In~\eqref{eq:proj_disc2}, $\mathbf{u}$ is a manufactured solution that coerces $\mathbf{u}_h$ to $\mathbf{u}$ and permits the \editor{%
truncation error and discretization error to be measured.  
Let $\mathbf{u}_n(\mathbf{x}) =\sum_{\srcidx=1}^{\nbasis} u_{\srcidx}^n \boldsymbol{\phi}_{\srcidx}(\mathbf{x})$ denote the basis-function representation of $\mathbf{u}$.  The truncation error is 
\begin{align}
\tau_{h_\testidx}(\mathbf{u})=a(\mathbf{u}_n,\boldsymbol{\phi}_{\testidx})-\left(\mathbf{f}, \boldsymbol{\phi}_{\testidx}\right)=a(\mathbf{u}_n,\boldsymbol{\phi}_{\testidx})-a(\mathbf{u},\boldsymbol{\phi}_{\testidx}).
\label{eq:efie_trunc}
\end{align}
As an alternative to the continuous form of the discretization error $\mathbf{e}_h=\mathbf{u}_h-\mathbf{u}$, we can measure the discrete form of the discretization error
\begin{align}
\mathbf{e}^h=\mathbf{u}^h-\mathbf{u}^n.
\label{eq:efie_disc}
\end{align}
}%

However, the integrals in~\eqref{eq:proj_disc2} cannot be accurately computed due to the presence of the Green's function, which yields a singularity when $\mathbf{x}=\mathbf{x}'$.  In~\cite{freno_em_mms_2020}, this challenge is mitigated  by manufacturing the Green's function as well.  As a result, the matrix becomes practically singular, admitting infinite solutions \editor{$\mathbf{u}^h$}.  Therefore, \editor{$\mathbf{u}^h$} is chosen by selecting the closest choice to \editor{$\mathbf{u}^n$} that satisfies~\eqref{eq:proj_disc2}.  The implications for the truncation and discretization errors as code-verification metrics for such singular systems of equations are discussed in the next section.
\section{Singular Systems of Equations}
\label{sec:singular}

Given their similar properties, in~\cite{warnick_2008}, expectations for the RWG basis functions are based on those for rooftop basis functions.  Likewise, in this section, we consider the one-dimensional analogy with piecewise linear basis functions $\phi_\srcidx(x)$.  We restrict the solution to real numbers and set $S=[a,\,b]$ in~\eqref{eq:a}.  We consider the extreme case of $G(x,x')=1$, which yields a singular system of equations.  As a result,~\eqref{eq:proj_disc2} becomes: find $u_h\in H_0^1$, such that
\editor{
\begin{align}
a(u_h,\phi_i) = a(u,\phi_i)
\label{eq:1d_var}
\end{align}
for $i=1,\hdots,n$, and 
\begin{align}
a(u,\phi_\testidx) = \alpha \int_a^b \phi_\testidx(x) \int_a^b u(x') dx' dx + \beta \int_a^b \frac{d\phi_\testidx}{dx}(x) \int_a^b \frac{du}{dx'}(x') dx' dx.
\label{eq:a_long}
\end{align}
%
Since $\int_a^b \frac{d\phi_\testidx}{dx}(x) dx=0$,~\eqref{eq:a_long} is reduced to
%
%
\begin{align*}
a(u,\phi_\testidx) = \alpha \int_a^b \phi_\testidx(x) \int_a^b u(x') dx' dx.
\end{align*}
}%

For a uniform discretization, after dividing both sides by $\alpha h$, \editor{where $h=(b-a)/N$ and $N=n+1$ is the number of elements}, the system of equations for~\eqref{eq:1d_var} is
\begin{align}
h\left[
\begin{matrix}
1      & \cdots & 1      \\
\vdots & \ddots & \vdots \\
1      & \cdots & 1      \\
\end{matrix}
\right]
\left\{\begin{matrix}u_1^h \\ \vdots \\ u_n^h\end{matrix}\right\}
=
\int_a^b u(x)dx
\left\{\begin{matrix}1 \\ \vdots \\ 1\end{matrix}\right\},
\label{eq:fe_mat}
\end{align}
which can be written as $\boldsymbol{\mathcal{L}}_h(u_h) = \mathbf{A}\mathbf{u}^h-\mathbf{b} = \mathbf{0}$, where $\mathbf{A}=h\mathbf{1}_{n\times n}$ and $\mathbf{b}=
\left(\int_a^b u(x)dx\right) \mathbf{1}_{n\times 1}$.  \editor{For the truncation error, the division by $h$ ensures that
\begin{align*}
\left.\frac{\partial \mathcal{L}_\testidx}{\partial v}\right|_u\!\cdot e_h=\frac{1}{\alpha h}a(e_h,\phi_\testidx)=\frac{1}{h} \int_a^b \phi_\testidx(x) \int_a^b e_h(x') dx' dx = \int_a^b e_h(x) dx,
\end{align*}
such that the truncation and discretization errors are of the same order in~\eqref{eq:linearization}.}


\subsection{Original System} 
\editor{%
Accounting for the boundary conditions $u(a)=u(b)=0$, the trapezoidal-rule approximation of the integral in~\eqref{eq:fe_mat} is
\begin{align}
\int_a^b u(x)dx = h\sum_{\srcidx=1}^\nbasis u_\srcidx^n + \mathcal{O}(h^2),
\label{eq:trap}
\end{align}
}%
such that $\mathbf{b}=\mathbf{A}\editor{\mathbf{u}^n}+\mathcal{O}(h^2)$ and~\eqref{eq:fe_mat} becomes
\begin{align}
\boldsymbol{\mathcal{L}}_h(u_h) = 
\mathbf{A}\left(\mathbf{u}^h-\editor{\mathbf{u}^n}\right)
{}+{}
\mathcal{O}(h^2)=\mathbf{0}.
\label{eq:fe_trunc}
\end{align}
From~\eqref{eq:trunc1} and~\eqref{eq:fe_trunc},
\begin{align}
\boldsymbol{\tau}_h(u) = \boldsymbol{\mathcal{L}}_h(u) = \editor{\tilde{\boldsymbol{\tau}}_h(u) +{}} \mathcal{O}(h^2),
\label{eq:fe_trunc2}
\end{align}
where 
\begin{align*}
\tilde{\boldsymbol{\tau}}_h(u)=\mathbf{A}\left(\mathbf{u}^n-\mathbf{u}^n\right)=\mathbf{0},
\end{align*}
such that the truncation error \reviewerTwo{$\boldsymbol{\tau}_h(u)$} is $\mathcal{O}(h^2)$, as expected.  Because $\mathbf{A}$ is singular, a unique solution cannot be computed and a meaningful discretization error cannot be obtained.

If an error $\editor{\delta}\ne 0$ \reviewerTwo{that is $\mathcal{O}(h^r)$, where $r\ge 0$,} is introduced in $A_{i,j}$ in~\eqref{eq:fe_mat}, such that it becomes $(1+\editor{\delta})A_{i,j}$, \editor{$\tilde{\boldsymbol{\tau}}_h(u)$ in the truncation error~\eqref{eq:fe_trunc2} becomes
%
%
\begin{align}
\tilde{\boldsymbol{\tau}}_h(u) = \delta u_j^n h\left\{\begin{array}{@{} l @{}} \mathbf{0}_{i-1\times 1} \\  \mathbf{1}_{1\times 1}\\  \mathbf{0}_{n-i\times 1} \end{array}\right\}.
\label{eq:tau_tilde}
\end{align}
}%
If $j$ is fixed, as $n$ is increased, $x_j$ will approach $a$ and $\editor{u_j^n}=u(x_j)$ will be \editor{asymptotically} proportional to the leading term of the Taylor series expansion at \editor{$x=a+jh$} about $x=a$.  If that term is $\mathcal{O}(h^q)$, \reviewerTwo{$\boldsymbol{\tau}_h(u)$}~\eqref{eq:fe_trunc2} will \editor{be} $\mathcal{O}(h^{\min\{\reviewerTwo{q+r+1},\,2\}})$.  If the error is introduced in a column corresponding to a fixed $x$ location (e.g., \editor{$j=N/2$}), \reviewerTwo{$\boldsymbol{\tau}_h(u)$} will \editor{be} $\mathcal{O}(h^\reviewerTwo{r+1})$.  Therefore, if \reviewerTwo{$j$ is fixed and $q+r\ge 1$ or the error is spatially fixed and $r\ge 1$, the error will go undetected.}

\subsection{Reduced System} 
To perform code verification for a singular system of equations, such as that in~\eqref{eq:fe_mat}, an approach is presented in~\cite{freno_em_mms_2020}, which uses a QR factorization of the transpose of the matrix to select the $\editor{\mathbf{u}^h}$ closest to $\editor{\mathbf{u}^n}$ that satisfies the system of equations. 
To simplify the factorization, we divide both sides of~\eqref{eq:fe_mat} by $h$ once more:
\begin{align}
\left[
\begin{matrix}
1      & \cdots & 1      \\
\vdots & \ddots & \vdots \\
1      & \cdots & 1      \\
\end{matrix}
\right]
\left\{\begin{matrix}u_1^h \\ \vdots \\ u_n^h\end{matrix}\right\}
=
\frac{1}{h}\int_a^b u(x)dx
\left\{\begin{matrix}1 \\ \vdots \\ 1\end{matrix}\right\},
\label{eq:fe_mat2}
\end{align} 
which can be written as $\mathbf{A}\mathbf{u}^h=\mathbf{b}$, where $\mathbf{A}=\mathbf{1}_{n\times n}$ and $\mathbf{b}=\frac{1}{h}\left(\int_a^b u(x) dx\right) \mathbf{1}_{n\times 1}$.  Performing a pivoted QR factorization of $\mathbf{A}^T$, such that
\begin{align*}
\mathbf{A}^T\mathbf{P} = \left[\mathbf{Q}_1,\, \mathbf{Q}_2\right]\left[\begin{matrix}\mathbf{R}_1 \\ \mathbf{0}\end{matrix}\right] = \mathbf{Q}_1\mathbf{R}_1,
\end{align*}
the optimal solution is
\begin{align}
\editor{\mathbf{u}^h} = \mathbf{Q}_1\mathbf{u}' + \mathbf{Q}_2\mathbf{Q}_2^T\editor{\mathbf{u}^n},
\label{eq:opta}
\end{align}
where $\mathbf{u}'$ is the solution to~\eqref{eq:fe_mat2} when the size is reduced to its rank:
\begin{align*}
\mathbf{u}' = \left(\mathbf{R}_1^T\right)^\dagger\mathbf{P}^T\mathbf{b},  
\end{align*}
where $(\cdot)^\dagger$ denotes the Moore--Penrose pseudoinverse.  Noting that $\mathbf{Q}\mathbf{Q}^T =\mathbf{Q}_1 \mathbf{Q}_1^T + \mathbf{Q}_2\mathbf{Q}_2^T=\mathbf{I}_{n\times n}$,~\eqref{eq:opta} can be written as
\begin{align}
\editor{\mathbf{u}^h} = \mathbf{Q}_1\mathbf{u}' - \mathbf{Q}_1\mathbf{Q}_1^T\editor{\mathbf{u}^n}+ \editor{\mathbf{u}^n}.
\label{eq:opt}
\end{align}

\subsubsection{Without an Error} 
The factorization is $\mathbf{P}=\mathbf{I}_{n\times n}$, $\mathbf{Q}_1 = n^{-1/2}\mathbf{1}_{n\times 1}$, $\mathbf{R}_1 = n^{1/2}\mathbf{1}_{1\times n}$, and $\left(\mathbf{R}_1^T\right)^\dagger=n^{-3/2}\mathbf{1}_{1\times n}$.
The first term of~\eqref{eq:opt} is 
\begin{align}
\mathbf{Q}_1\mathbf{u}'=\mathbf{Q}_1\left(\mathbf{R}_1^T\right)^\dagger\mathbf{P}^T\mathbf{b} = \frac{1}{n^2}\mathbf{1}_{n\times n} \mathbf{b} = \frac{1}{nh}\left(\int_a^b u(x)dx\right) \mathbf{1}_{n\times 1}.
\label{eq:Q1u1}
\end{align}
Inserting~\eqref{eq:trap} into~\eqref{eq:Q1u1} yields
\begin{align}
\mathbf{Q}_1\mathbf{u}' &{}=%
\frac{1}{n}\mathbf{1}_{n\times n}\editor{\mathbf{u}^n} + \mathcal{O}(h^2).
\label{eq:Q1up}
\end{align}
The second term of~\eqref{eq:opt} is 
\begin{align}
\mathbf{Q}_1\mathbf{Q}_1^T\editor{\mathbf{u}^n} =\frac{1}{n}\mathbf{1}_{n\times n}
\editor{\mathbf{u}^n}.
\label{eq:Q2Q2u}
\end{align}

Inserting~\eqref{eq:Q1up} and~\eqref{eq:Q2Q2u} into~\eqref{eq:opt}, the discretization error~\eqref{eq:efie_disc} is
\editor{%
\begin{align}
\mathbf{e}^h &{}= \tilde{\mathbf{e}}^h + \mathcal{O}(h^2),
\label{eq:disc_good}
\end{align}
where
\begin{align*}
\tilde{\mathbf{e}}^h =
\left(%
\frac{1}{n}\mathbf{1}_{n\times n} - \frac{1}{n}\mathbf{1}_{n\times n}
{}+{}
\mathbf{I}_{n\times n}
{}-{}
\mathbf{I}_{n\times n}
\right)\mathbf{u}^n=\mathbf{0}.
\end{align*}
}%
The discretization error \reviewerTwo{$\mathbf{e}^h$} is $\mathcal{O}(h^2)$, as expected.

\subsubsection{With an Error} 
We now consider the case when an error $\editor{\delta}\ne 0$ is introduced in $A_{i,j}$ in~\eqref{eq:fe_mat2}, such that it becomes $(1+\editor{\delta})A_{i,j}$.  If $i=1$, the factorization is $\mathbf{P}=\mathbf{I}_{n\times n}$,
\begin{align*}
\mathbf{Q}_1 &{}= \frac{1}{\gamma}\left[\begin{array}{@{} r @{} l r @{} l @{}}
               &\mathbf{1}_{j-1\times 1} &  \ope\nmo^{-1/2           }&\mathbf{1}_{j-1\times 1} \\
\ope&\mathbf{1}_{1\times 1}   &                -\nmo^{ 1/2\phantom{-}}&\mathbf{1}_{1  \times 1} \\
               &\mathbf{1}_{n-j\times 1} &  \ope\nmo^{-1/2           }&\mathbf{1}_{n-j\times 1}
\end{array}\right],
\\[1.5em]
\mathbf{R}_1 &{}= \frac{1}{\gamma}\left[\begin{array}{@{} r @{} l r @{} l @{}}
\gamma^2&\mathbf{1}_{1\times 1} & \npe&\mathbf{1}_{1\times n-1} \\
                                  &\mathbf{0}_{1\times 1} & \editor{\delta}\nmo^{1/2}    &\mathbf{1}_{1\times n-1}
\end{array}\right],
\\[1.5em] 
\left(\mathbf{R}_1^T\right)^\dagger &{}= \frac{1}{\editor{\delta}\gamma}\left[\begin{array}{@{} r @{} l r @{} l @{}}
\editor{\delta}&\mathbf{1}_{1\times 1} & &\mathbf{0}_{1\times n-1} \\
-\npe\nmo^{-1/2}&\mathbf{1}_{1\times 1} &   \gamma^2\nmo^{-3/2}&\mathbf{1}_{1\times n-1}
\end{array}\right],
\end{align*}
where $\nmo=n-1$, $\ope=1+\editor{\delta}$, $\npe=n+\editor{\delta}$, and $\gamma=\sqrt{\nmo+\ope^2}$.
If $i\ne 1$, the factorization is $\mathbf{P}=\mathbf{I}_{n\times n}$,
\begin{align*}
\mathbf{Q}_1 &{}= n^{-1/2}\left[\begin{array}{@{} r @{} l r @{} l @{}}
&\mathbf{1}_{j-1\times 1} & -\nmo^{-1/2}          &\mathbf{1}_{j-1\times 1} \\
&\mathbf{1}_{1\times 1}   &  \nmo^{1/2\phantom{-}}&\mathbf{1}_{1\times 1}\\
&\mathbf{1}_{n-j\times 1} & -\nmo^{-1/2}          &\mathbf{1}_{n-j\times 1}
\end{array}\right],
\\[1.5em]
\mathbf{R}_1 &{}= n^{-1/2}\left[\begin{array}{@{} r @{} l r @{} l r @{} l @{}}
n&\mathbf{1}_{1\times i-1} &  \npe&\mathbf{1}_{1\times 1} & n&\mathbf{1}_{1\times n-i} \\
 &\mathbf{0}_{1\times i-1} & \editor{\delta}\nmo^{1/2}    &\mathbf{1}_{1\times 1} & &\mathbf{0}_{1\times n-i}
\end{array}\right],
\\[1.5em]
\left(\mathbf{R}_1^T\right)^\dagger &{}=\frac{n^{-1/2}\nmo^{-3/2}}{\editor{\delta}}\left[\begin{array}{@{} r @{} l r @{} l r @{} l @{}}
\editor{\delta}\nmo^{1/2}&\mathbf{1}_{1\times i-1} &  &\mathbf{0}_{1\times 1} & \editor{\delta}\nmo^{1/2}&\mathbf{1}_{1\times n-i} \\
-\npe&\mathbf{1}_{1\times i-1} & n\nmo &\mathbf{1}_{1\times 1} & -\npe&\mathbf{1}_{1\times n-i}
\end{array}\right].
\end{align*}

The first term of~\eqref{eq:opt} is
\begin{align}
\mathbf{Q}_1\mathbf{u}'=\mathbf{Q}_1\left(\mathbf{R}_1^T\right)^\dagger\mathbf{P}^T\mathbf{b} &{}= 
\frac{1}{\editor{\delta}\nmo^2}
\left[\begin{array}{@{} r @{} l r @{} l r @{} l @{}}
\ope & \mathbf{1}_{j-1\times i-1} & -\nmo\phantom{^1}&\mathbf{1}_{j-1\times 1}  & \ope& \mathbf{1}_{j-1\times n-i} \\
-\nmo          & \mathbf{1}_{1\times i-1}   &  \nmo^2          &\mathbf{1}_{1\times 1}  & -\nmo         & \mathbf{1}_{1\times n-i}   \\
\ope & \mathbf{1}_{n-j\times i-1} & -\nmo\phantom{^1}&\mathbf{1}_{n-j\times 1}  & \ope& \mathbf{1}_{n-j\times n-i}
\end{array}\right]\mathbf{b} \nonumber \\[1em]
&{}= 
\frac{1}{\nmo h}\int_a^b u(x)dx
\left\{\begin{array}{@{} l @{}}
\mathbf{1}_{j-1\times 1} \\ \mathbf{0}_{1\times 1} \\ \mathbf{1}_{n-j\times 1}
\end{array}\right\}.
\label{eq:Q1u1_bad}
\end{align}
Inserting~\eqref{eq:trap} into~\eqref{eq:Q1u1_bad} yields
\begin{align}
\mathbf{Q}_1\mathbf{u}'=\frac{1}{\nmo}
\left[\begin{array}{@{} l @{}}
\mathbf{1}_{j-1\times n} \\ \mathbf{0}_{1 \times n} \\ \mathbf{1}_{n-j\times n}
\end{array}\right]
\editor{\mathbf{u}^n} + \mathcal{O}(h^2).
\label{eq:Q1up_bad}
\end{align}
The second term of~\eqref{eq:opt} is 
\begin{align}
\mathbf{Q}_1\mathbf{Q}_1^T\editor{\mathbf{u}^n} &{}=
\left(\left[\begin{array}{@{} l l l @{}}
\mathbf{0}_{j-1\times j-1} & \mathbf{0}_{j-1\times 1} & \mathbf{0}_{j-1\times n-j} \\
\mathbf{0}_{  1\times j-1} & \mathbf{1}_{  1\times 1} & \mathbf{0}_{1  \times n-j} \\
\mathbf{0}_{n-j\times j-1} & \mathbf{0}_{n-j\times 1} & \mathbf{0}_{n-j\times n-j}
\end{array}\right]
{}+{}
\frac{1}{\nmo}
\left[\begin{array}{@{} l ll @{}}
\mathbf{1}_{j-1\times j-1} & \mathbf{0}_{j-1\times 1} & \mathbf{1}_{j-1\times n-j} \\
\mathbf{0}_{  1\times j-1} & \mathbf{0}_{  1\times 1} & \mathbf{0}_{  1\times n-j} \\
\mathbf{1}_{n-j\times j-1} & \mathbf{0}_{n-j\times 1} & \mathbf{1}_{n-j\times n-j}
\end{array}\right]
\right)\editor{\mathbf{u}^n}.
\label{eq:Q2Q2u_bad}
\end{align}

Inserting~\eqref{eq:Q1up_bad} and~\eqref{eq:Q2Q2u_bad} into~\eqref{eq:opt}, $\tilde{\mathbf{e}}^h$ in the discretization error \reviewerTwo{$\mathbf{e}^h = \tilde{\mathbf{e}}^h + \mathcal{O}(h^2)$}~\eqref{eq:disc_good} is
\begin{align}
\tilde{\mathbf{e}}^h &{}= 
\left(-%
\left[\begin{array}{@{} l l l @{}}
\mathbf{0}_{j-1\times j-1} & \mathbf{0}_{j-1\times 1} & \mathbf{0}_{j-1\times n-j} \\
\mathbf{0}_{  1\times j-1} & \mathbf{1}_{  1\times 1} & \mathbf{0}_{1  \times n-j} \\
\mathbf{0}_{n-j\times j-1} & \mathbf{0}_{n-j\times 1} & \mathbf{0}_{n-j\times n-j}
\end{array}\right]
{}+{}
\frac{1}{\nmo}
\left[\begin{array}{@{} l l l @{}}
\mathbf{0}_{j-1\times j-1} & \mathbf{1}_{j-1\times 1} & \mathbf{0}_{j-1\times n-j} \\
\mathbf{0}_{  1\times j-1} & \mathbf{0}_{  1\times 1} & \mathbf{0}_{  1\times n-j} \\
\mathbf{0}_{n-j\times j-1} & \mathbf{1}_{n-j\times 1} & \mathbf{0}_{n-j\times n-j}
\end{array}\right]
\right)\editor{\mathbf{u}^n}
\nonumber \\[1em]
&{}= 
u_j^n\left(-\left\{\begin{array}{@{} l @{}}
\mathbf{0}_{j-1\times 1} \\ \mathbf{1}_{1\times 1} \\ \mathbf{0}_{n-j\times 1}
\end{array}\right\} 
{}+{}
\frac{h}{b-a}
\left\{\begin{array}{@{} l @{}}
\mathbf{1}_{j-1\times 1} \\ \mathbf{0}_{1\times 1} \\ \mathbf{1}_{n-j\times 1}
\end{array}\right\}
\right).
\label{eq:eh_tilde}
\end{align}
If $j$ is fixed, as $n$ is increased, $x_j$ will approach $a$ and $\editor{u_j^n}=u(x_j)$ will be \editor{asymptotically} proportional to the leading term of the Taylor series expansion at \editor{$x=a+jh$} about $x=a$.  If that term is $\mathcal{O}(h^q)$,~\eqref{eq:disc_good} will \editor{be} $\mathcal{O}(h^{\min\{q,\,2\}})$.  If the error is introduced in a column corresponding to a fixed $x$ location,~\eqref{eq:disc_good} will \editor{be} $\mathcal{O}(h^0)$ \editor{if $u(x_j)\ne 0$}.  Therefore, if $q<2$ or the error is spatially fixed, the error will be detected.  \editor{Note that, while the column index of the error appears in~\eqref{eq:eh_tilde}, $\delta$ does not.}

\subsection{\reviewerTwo{Summary}}

\reviewerTwo{%
In the absence of coding errors, the truncation error computed from the singular equations and the discretization error computed from the optimal solution are both $\mathcal{O}(h^2)$.  However, the discretization error is more effective at detecting coding errors than the truncation error.  

For coding errors with an index-fixed column, the discretization error can detect errors if $q<2$, whereas the truncation error can only detect errors if $q< 1 - r$.  For coding errors with a physically fixed column, the discretization error can detect any errors, whereas the truncation error can only detect errors if $r<1$.  
}
\section{Numerical Examples}

We demonstrate the observations of Section~\ref{sec:singular} with numerical examples in this section for the simplified one-dimensional analogy to the EFIE~\eqref{eq:fe_mat}, as well as for the EFIE~\eqref{eq:proj_disc2}.  

\subsection{Simplified Analogy} 

\begin{table}
\centering
\editor{%
\begin{tabularx}{0.7\textwidth}{c c c >{\centering}X >{\centering}X >{\centering}X >{\centering\arraybackslash}X }
\toprule
     &                         &     & \multicolumn{4}{c}{Location of $\delta$, $(i,j)$} \\ \cmidrule(){4-7} 
Case & $u(x)$                  & $q$ &  a  &  b      &   c       & d                     \\ \midrule
1    & $\sin\pi x\phantom{^1}$ & 1   & --- & $(1,2)$ & $(2,N/4)$ & $(N/2,N/2)$           \\
2    & $\sin\pi x^2$           & 2   & --- & $(1,2)$ & $(2,N/4)$ & $(N/2,N/2)$           \\
\bottomrule
\end{tabularx}}
\caption{\editor{Simplified analogy: Cases.}}
\label{tab:cases}
\end{table}

\begin{figure}
\centering
\begin{subfigure}[t]{.49\textwidth}
\includegraphics[scale=.64,clip=true,trim=2.25in 0in 2.8in 0in]{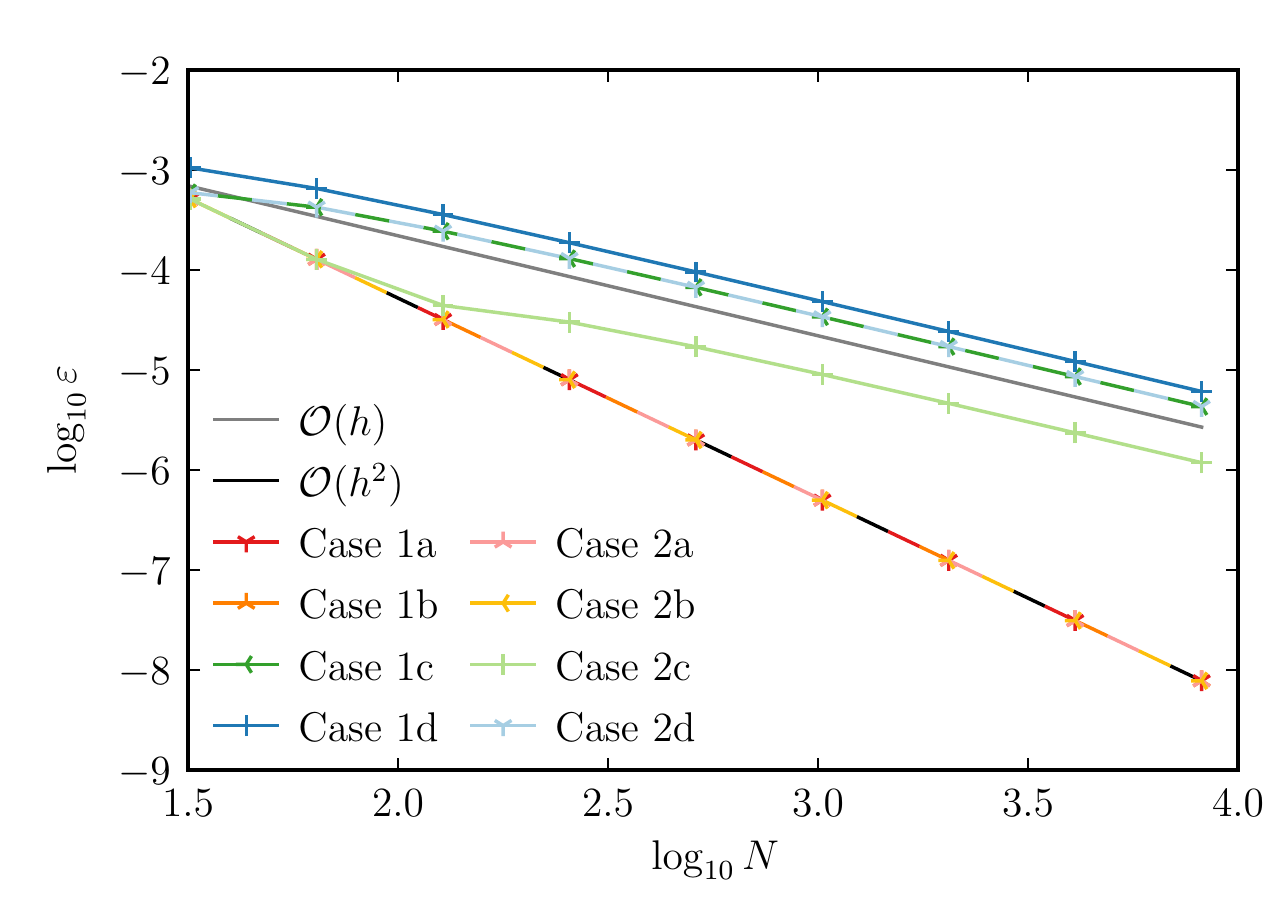}
\caption{$\varepsilon={\|\boldsymbol{\tau}_h(u)\|}_\infty$}
\label{fig:trunc_linf}
\end{subfigure}
\hspace{0.25em}
\begin{subfigure}[t]{.49\textwidth}
\includegraphics[scale=.64,clip=true,trim=2.25in 0in 2.8in 0in]{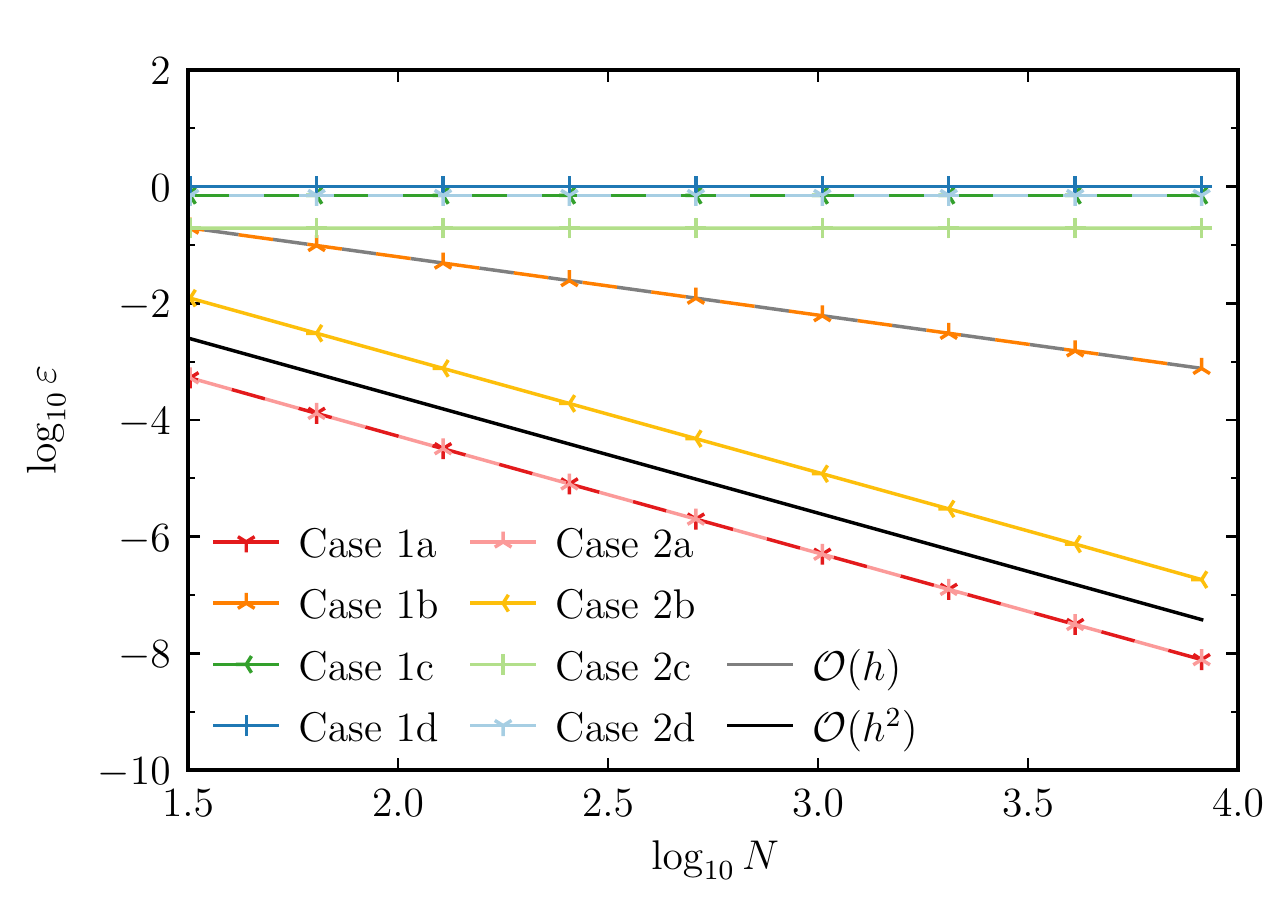}
\caption{$\varepsilon={\|\mathbf{e}^h\|}_\infty$}
\label{fig:disc_linf}
\end{subfigure}
\\
\begin{subfigure}[t]{.49\textwidth}
\includegraphics[scale=.64,clip=true,trim=2.25in 0in 2.8in 0in]{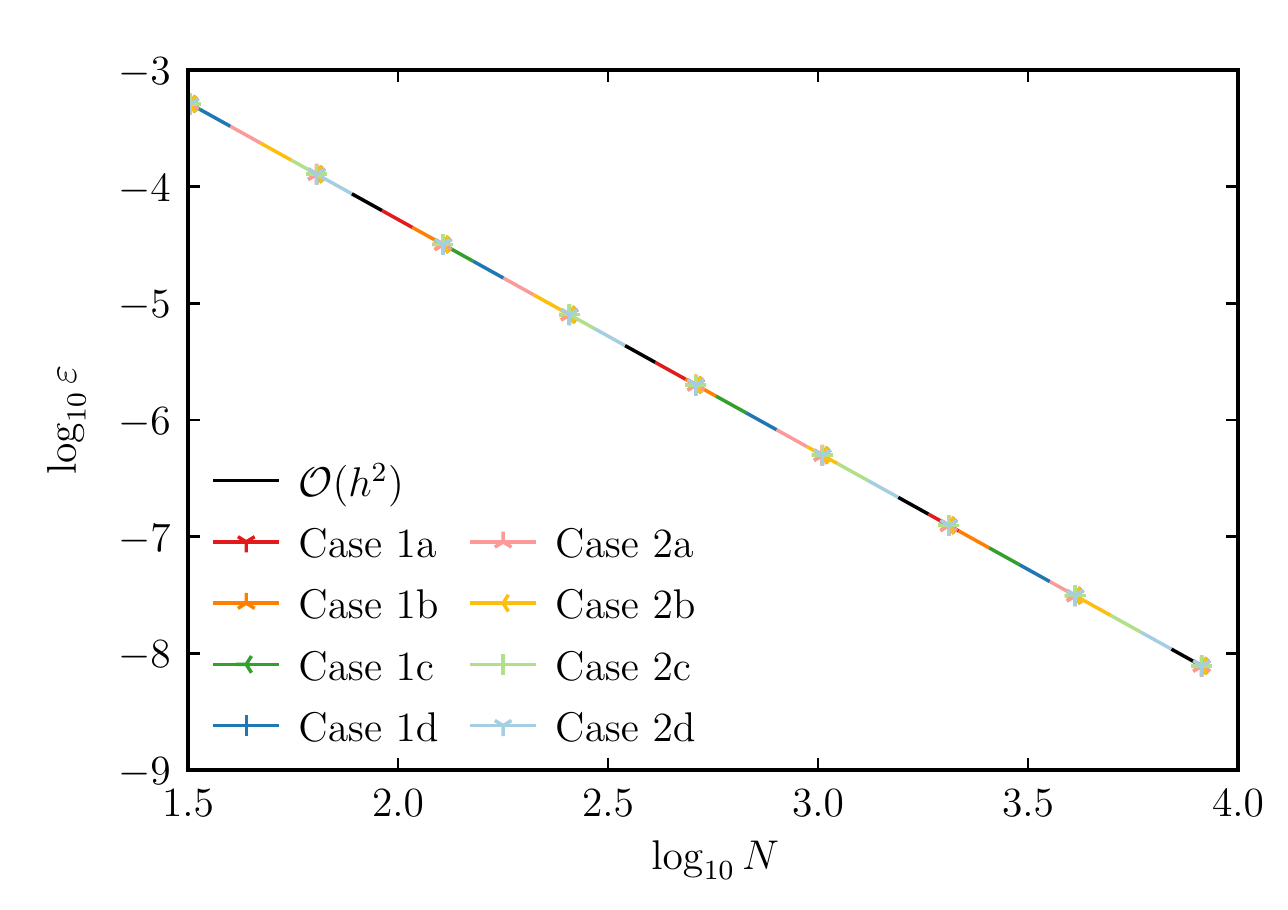}
\caption{$\varepsilon={\|\tilde{\boldsymbol{\tau}}_h(u)-\boldsymbol{\tau}_h(u)\|}_\infty$}
\label{fig:err_trunc_linf}
\end{subfigure}
\hspace{0.25em}
\begin{subfigure}[t]{.49\textwidth}
\includegraphics[scale=.64,clip=true,trim=2.25in 0in 2.8in 0in]{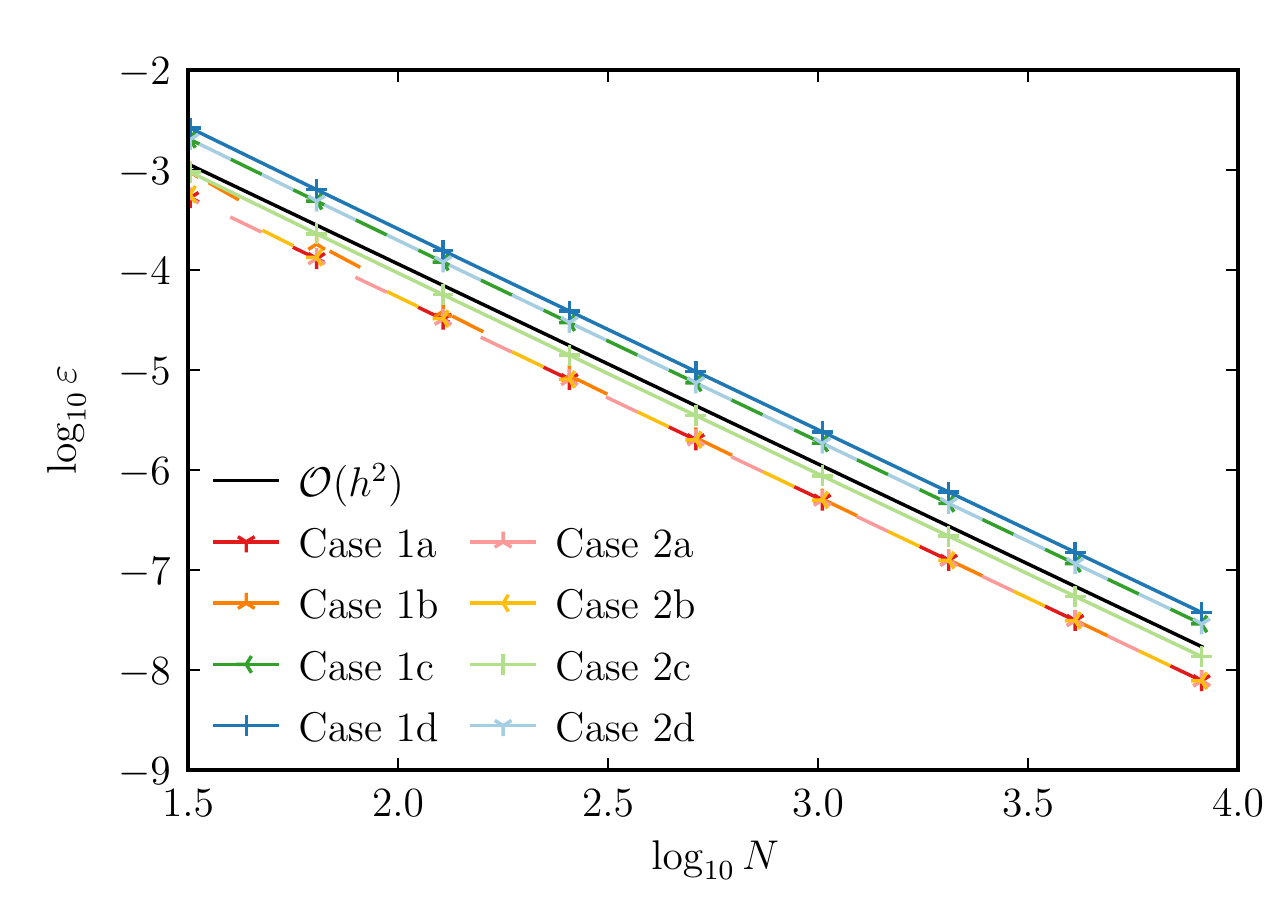}
\caption{$\varepsilon={\|\tilde{\mathbf{e}}^h -\mathbf{e}^h \|}_\infty$}
\label{fig:err_disc_linf}
\end{subfigure}
\caption{\editor{Simplified analogy: Convergence comparison of error metrics $\varepsilon$.}}
\vskip-\dp\strutbox
\label{fig:error_metrics}
\end{figure}

For the simplified analogy, we consider the cases in Table~\ref{tab:cases} for the domain $[a,\,b]=[0,\,1]$.  Cases 1a and 2a have no errors, Cases 1b and 2b have an error of $\delta=1/20$ at fixed row and column indices, Cases 1c and 2c have the error at a fixed row index and spatially fixed column, and Cases 1d and 2d have the error at spatially fixed rows and columns.  The leading term of the Taylor series expansion at $x = a+h$ about $x=a$ is $\mathcal{O}(h^q)$, where $q=1$ for Case 1 and $q=2$ for Case 2.

Figure~\ref{fig:error_metrics} provides a comparison of the error metrics $\boldsymbol{\tau}_h(u)$ from~\eqref{eq:fe_trunc2}, $\mathbf{e}^h$ from~\eqref{eq:disc_good}, $\tilde{\boldsymbol{\tau}}_h(u)-\boldsymbol{\tau}_h(u)$ from~\eqref{eq:fe_trunc2}, and $\tilde{\mathbf{e}}^h -\mathbf{e}^h$ from~\eqref{eq:disc_good}.
$\boldsymbol{\tau}_h(u)$ is $\mathcal{O}(h^2)$ not only for Cases 1a and 2a, but also for Cases 1b and 2b, such that the fixed-index errors go undetected.  Cases 1c, 1d, 2c, and 2d are $\mathcal{O}(h)$, such that the spatially fixed errors are detected, but consistency is falsely implied.
$\mathbf{e}^h$ is $\mathcal{O}(h^2)$ for Cases 1a and 2a, as well as Case 2b since the error introduced in Case 2b is $\mathcal{O}(h^2)$ since $q=2$.  Case 1b is $\mathcal{O}(h)$ and Cases 1c, 1d, 2c and 2d are $\mathcal{O}(h^0)$, such that those errors are detected.
Finally, $\tilde{\boldsymbol{\tau}}_h(u)-\boldsymbol{\tau}_h(u)$ verifies the expression in~\eqref{eq:tau_tilde} and $\tilde{\mathbf{e}}^h -\mathbf{e}^h$ verifies the expression in~\eqref{eq:eh_tilde}, as only the $\mathcal{O}(h^2)$ terms remain.

\subsection{EFIE} 

\begin{figure}
\centering
\begin{subfigure}[t]{.49\textwidth}
\includegraphics[scale=.64,clip=true,trim=2.25in 0in 2.8in 0in]{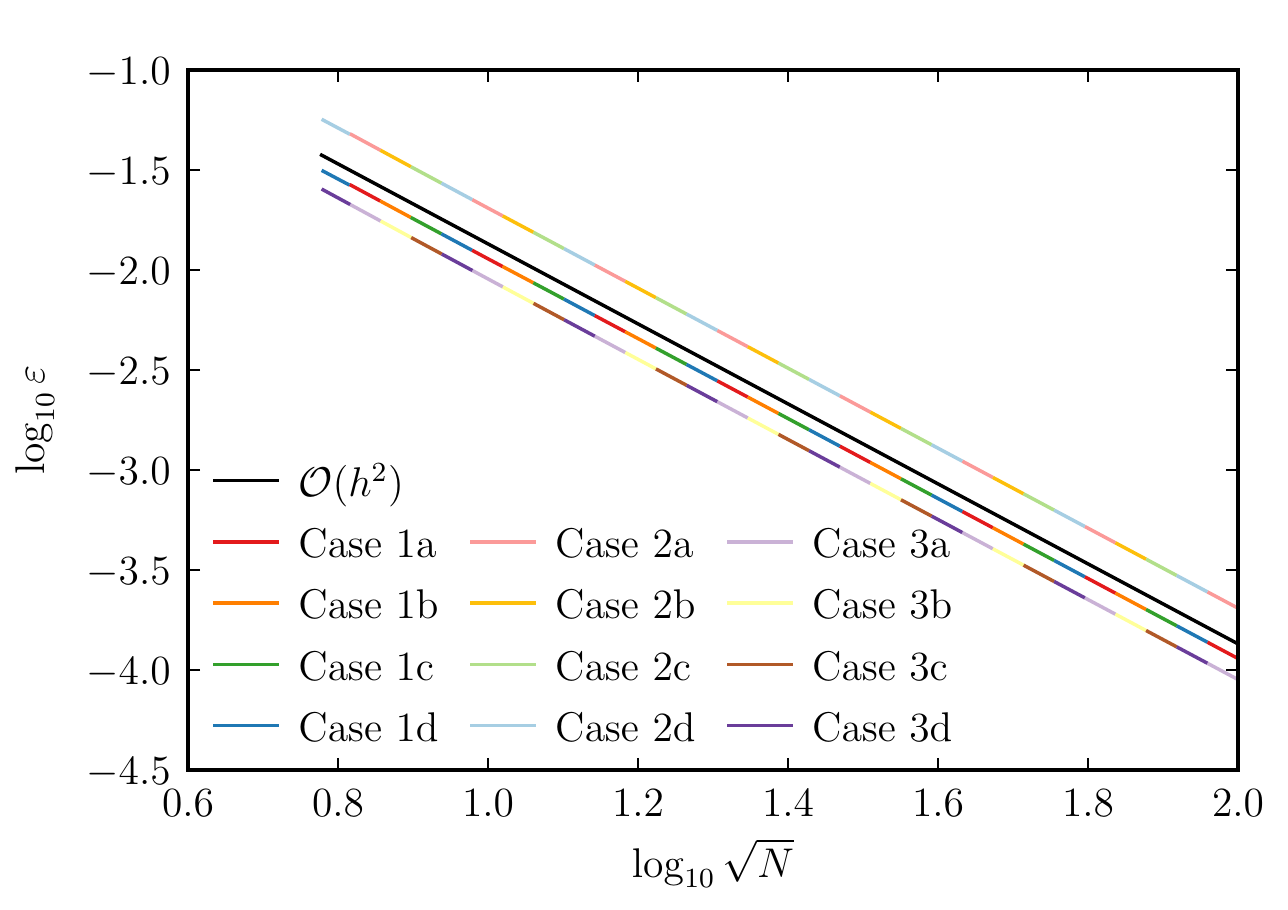}
\caption{$\varepsilon={\|\boldsymbol{\tau}_h(\mathbf{u})\|}_\infty$}
\label{fig:trunc_linf}
\end{subfigure}
\hspace{0.25em}
\begin{subfigure}[t]{.49\textwidth}
\includegraphics[scale=.64,clip=true,trim=2.25in 0in 2.8in 0in]{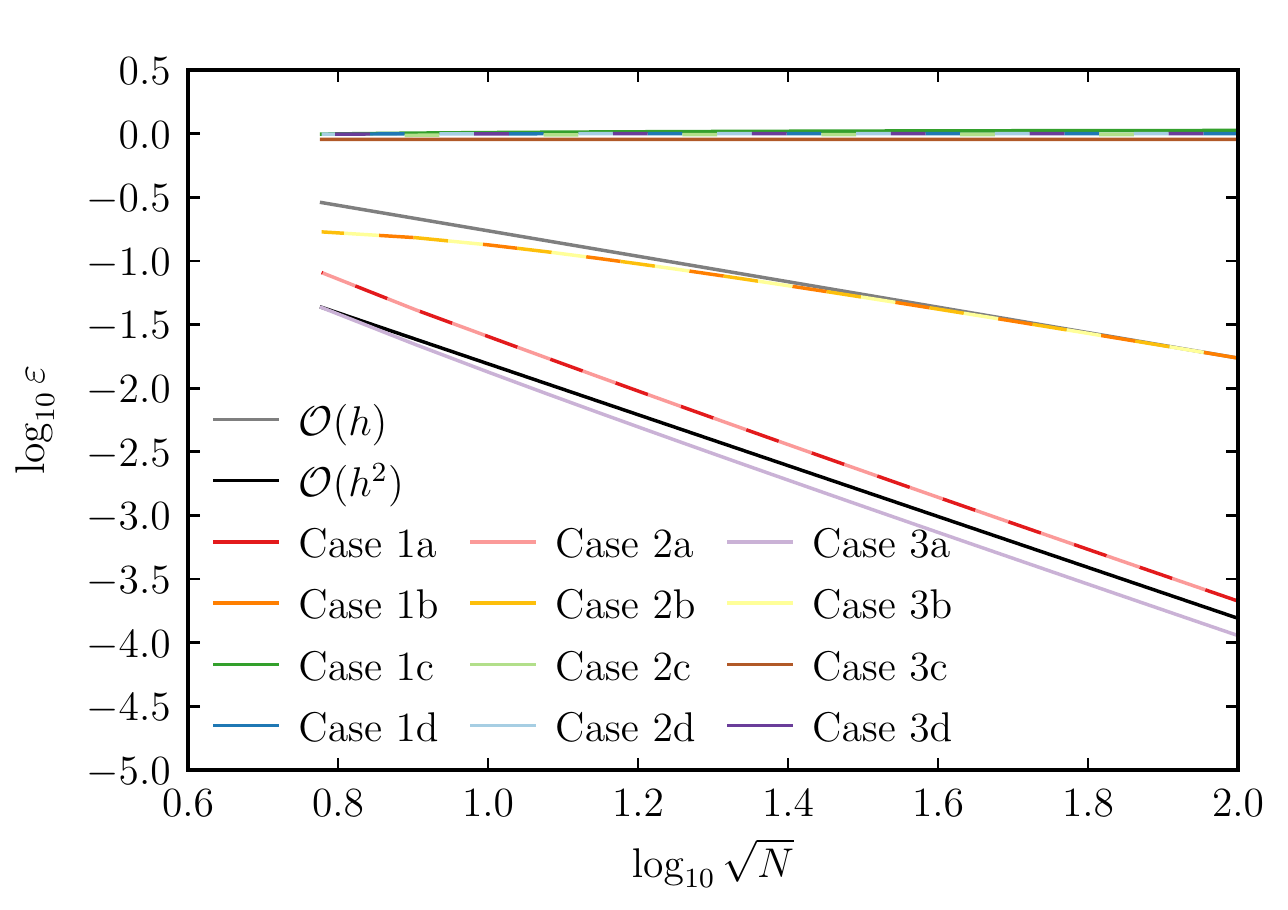}
\caption{$\varepsilon={\|\mathbf{e}^h\|}_\infty$}
\label{fig:disc_linf}
\end{subfigure}
\caption{\editor{EFIE: Convergence comparison of error metrics $\varepsilon$.}}
\vskip-\dp\strutbox
\label{fig:efie_error_metrics}
\end{figure}

\begin{table}
\centering
\editor{%
\begin{tabularx}{0.7\textwidth}{c c >{\centering}X >{\centering}X >{\centering}X >{\centering\arraybackslash}X }
\toprule
     &                  & \multicolumn{4}{c}{Location of $\delta$, $(i,j)$} \\ \cmidrule(){3-6} 
Case & $(\alpha,\beta)$ &  a  &  b      &   c        & d                    \\ \midrule
1    & $(1,1)$          & --- & $(1,2)$ & $(2,f(N))$ & $(f(N),f(N))$        \\
2    & $(1,0)$          & --- & $(1,2)$ & $(2,f(N))$ & $(f(N),f(N))$        \\
3    & $(0,1)$          & --- & $(1,2)$ & $(2,f(N))$ & $(f(N),f(N))$        \\
\bottomrule
\end{tabularx}}
\caption{\editor{EFIE: Cases.}}
\label{tab:efie_cases}
\end{table}

For the EFIE, we consider the cases in Table~\ref{tab:efie_cases} for the domain $(x,y)\in[-1,\,1]\times[0,\,1]$.  For these cases, $\mathbf{u}(\mathbf{x})=\{\cos\left(\pi x/2\right)\cos\left(\pi y/4\right),\cos(\pi x/4)\sin(\pi y)\}$, and $G(\mathbf{x},\mathbf{x}')=1-\|\mathbf{x}-\mathbf{x}'\|_2^2/5$.  If $j$ is fixed, as $n$ is increased, $(x_j,y_j)$ will approach $(-1,0)$, and the leading term of the Taylor series expansion is $\mathcal{O}(h)$.  For the EFIE, $N$ denotes the number of triangles used to discretize the domain and $h=2/\sqrt{N}$.

Cases 1a, 2a, and 3a have no errors.  
Cases 1b, 2b, and 3b have an error of $\delta=1/100$ at fixed row and column indices.
For Cases 1c, 2c, and 3c, the row index of the error is fixed, and the column indices $f(N)$ are those associated with the solution values at $x=0$. For Cases 1d, 2d, and 3d, the row and column indices of the error are spatially fixed and denoted by $f(N)$.

Figure~\ref{fig:efie_error_metrics} provides a comparison of the error metrics $\boldsymbol{\tau}_h(\mathbf{u})$ from~\eqref{eq:efie_trunc} and $\mathbf{e}^h$ from~\eqref{eq:efie_disc}.
For the EFIE, the truncation error is computed by factoring out $h^2$, so that 
\begin{align*}
\left.\frac{\partial \mathcal{L}_\testidx}{\partial\mathbf{v}}\right|_\mathbf{u}\!\cdot\mathbf{e}_h=\frac{a(\mathbf{e}_h,\boldsymbol{\phi}_\testidx)}{h^2},
\end{align*}
and the truncation and discretization errors are of the same order in~\eqref{eq:linearization}.  For the cases with $\delta$, 
\begin{align}
\tau_{h_{\testidx}}(\mathbf{u}) = \frac{\delta u_j^n}{h^2}a(\boldsymbol{\phi}_{\srcidx},\boldsymbol{\phi}_{\testidx})+ \mathcal{O}(h^2).
\label{eq:efie_tau}
\end{align}
In~\eqref{eq:efie_tau}, $a(\boldsymbol{\phi}_{\srcidx},\boldsymbol{\phi}_{\testidx})$ is $\mathcal{O}(h^4)$ due to the four-dimensional integrals of the basis functions.
Therefore, $\boldsymbol{\tau}_h(\mathbf{u})$ is $\mathcal{O}(h^2)$ not only for the cases without error, but for all of the cases with error, such that none of the errors are detected.
$\mathbf{e}^h$ is $\mathcal{O}(h^2)$ for Cases 1a, 2a, and 3a.  Cases 1b, 2b, and 3b are $\mathcal{O}(h)$, and Cases 1c, 1d, 2c, 2d, 3c, and 3d are $\mathcal{O}(h^0)$, such that these errors are detected.

\section{Conclusions}

In this paper, we showed how, for a singular system of equations, computing the truncation error by inserting the exact solution into the discretized equations cannot detect certain orders of coding errors.  However, the discretization error from the optimal solution is a more effective metric.
\section*{Acknowledgments} 
\label{sec:acknowledgments}
This paper describes objective technical results and analysis. Any subjective views or opinions that might be expressed in the paper do not necessarily represent the views of the U.S. Department of Energy or the United States Government.
Sandia National Laboratories is a multimission laboratory managed and operated by National Technology and Engineering Solutions of Sandia, LLC, a wholly owned subsidiary of Honeywell International, Inc., for the U.S. Department of Energy's National Nuclear Security Administration under contract DE-NA-0003525.

%
\addcontentsline{toc}{section}{\refname}
\bibliographystyle{elsarticle-num}
\bibliography{../../quadrature_manuscript/quadrature.bib}

\end{document}